\documentclass[12pt]{amsart}
\usepackage{amsmath}

\usepackage{url}

 \makeatletter
\renewcommand*\subjclass[2][2010]{%
  \def\@subjclass{#2}%
  \@ifundefined{subjclassname@#1}{%
    \ClassWarning{\@classname}{Unknown edition (#1) of Mathematics
      Subject Classification; using '2010'.}%
  }{%
    \@xp\let\@xp\subjclassname\csname subjclassname@#1\endcsname  
}%
}

 \makeatother

\begin{document}

\title[Different classes of binary necklaces\ldots]{Different classes of binary 
necklaces and a combinatorial method for their enumerations}

\author{Romeo Me\v{s}trovi\'c}

\address{Maritime Faculty, University of Montenegro, Dobrota 36,
 85330 Kotor, Montenegro} \email{romeo@ac.me}

\maketitle

{\renewcommand{\thefootnote}{}

{\renewcommand{\thefootnote}{}
  \setcounter{footnote}{0}}

 \begin{abstract}
In this paper  we investigate  enumeration of 
some classes of $n$-character strings and   binary   necklaces.
Recall that binary necklaces are necklaces  in two colors with length 
$n$. 
We prove three results (Theorems 1, 1' and  2)
 concerning the numbers of three classes of $j$-character strings
(closely related to some classes of binary necklaces or Lyndon words). 
Using these results, we deduce Moreau's necklace-counting function for 
binary aperiodic necklaces  of length 
$k$ \cite{mo}  (Theorem 3), and we prove the binary case of MacMahon's  
formula from 1892 \cite{ma} (also called Witt's formula) for the number of necklaces (Theorem 4).  
Notice that we give 
proofs of Theorems 3 and 4   
without use of Burnside's lemma  and P\'olya enumeration theorem.
Namely, the methods used in our proofs of auxiliary and  
main results presented in Sections 3 and 4 are  
combinatorial in spirit  and they are based on counting method
and  some facts from elementary number theory. 
  \end{abstract}

\section{Introduction}


George P\'olya (1887--1985) discovered a powerful general method for 
enumerating the number of orbits of a group on particular configurations. 
This method became known as the P\'olya Enumeration Theorem, or PET, whose
proof follows directly from Burnside's lemma.  
     P\'olya's theorem can be used to enumerate several objectst under 
permutation groups. In particular, it can be used for enumeration of 
different classes of necklaces and bracelets.

In combinatorics, a $k$-{\it ary necklace} of length $n$ is an equivalence 
class 
of $n$-character string  over an alphabet of size $k$, taking all rotations 
are equivalent. It represents a structure with $n$ circularty connected 
beads of up to $k$ different colors. A necklace of length $n$ is 
{\it primitive} if its period is not a proper divisor of $n$. 
  
Technically, one may classify as an orbit of the action of the cyclic group 
of $n$-character strings, and a bracelet as an orbit of the dihedral grooup's
action. Namely, an $(n,k)$ - bracelet is an equivalence classs of words of 
length $n$ under rotation and reflection.  
This enables appplication of P\'olya enumeration theorem of necklaces 
and bracelets. An $(n,k)$ - necklace is an equivalence class of words of length 
$n$ over an alphabet of size $k$ under rotation. For example, if $k=2$  
and the alphabet is $\{0,1\}$, then the following sets are  examples of three 
 binary necklaces (i.e., those with $k=2$):
  $$
\{0101,1010\},
  $$ 
  $$
\{011011,110110,1011101\},
$$
and
 $$
\{0110101,1101010,1010101,0101011,1010110,0101101,1011010\},
    $$  
The basic enumeration problem is 
then (Necklace Enumeration): {\it For a given $n$ and $k$, how many 
$(n,k)$-necklaces are there}? Equivalently,  we are asking how many orbits 
 the  cyclic group $C_n$ has on the set of all words of  length $n$ over an 
alphabet of size $k$. We will denote this value by $a(n,k)$. Notice that in 
a group $G$ of symmetry transformations such that only translations 
$(a_i\to a_{i+s})$ are allowed, $G$ is a cyclic group $C_n$. This case appears 
in \cite{ri} in connection with counting necklaces made from $n$ beads 
of $k$  different kinds (translations merely rotate the necklace). It also 
arises in problems of coding and genetics \cite{ggw}. The special
case $n=12$, $k=2$ occurs in finding the number of distinct musical
chords (of 0, 1, $\cdots$, or  12 notes) when inversions and transpositions 
to other keys are equivalences (for related calculations see 
\cite[Section 6]{gir}).           
 
An aperiodic necklace of length $n$ is an equivalence class of size $n$, 
i.e., no distinct rotations of a necklace from such class are equal. 
According to Moreau's necklace-counting function (see \cite[p. 503]{lu}; also 
see \cite{ri}), 
there are 
     $$
M_k(n)=\frac{1}{n}\sum_{d\mid n}\mu(d)k^{n/d}\eqno(1)
  $$
different $k$-ary aperiodic necklaces  of length $n$, where $\mu$ is the 
M\"{o}bius function, where $\mu(1)=1$, $\mu(n)=(-1)^r$ if $n$ is a product 
of $r$ distinct primes, and $\mu(n)=0$ otherwise (see the sequence 
A001037 in   \cite{oe} concerning the sequence 
$\{M_2(n)\}_{n=1}^{\infty}$ which presents the number of {\it binary Lyndon 
words}). 
 The formula (1) is called  
 {\it MacMahon's} formula in the book by Graham et al. \cite[the formula 
(4.63), p. 141]{gr}.
Notice that this formula may be derived by a simple 
direct argument given in \cite{ggw}.

Each aperiodic necklace contains a single Lyndon word so that Lyndon words 
form representatives of aperiodic necklaces. Recall that in mathematics, in 
the areas of combinatorics and computer science, a {\it Lyndon word} is a 
nonempty string 
that is strictly smaller than lexicographic order than all of its rotations.
More precisely, a $k$-ary Lyndon word of length $n>0$ is an $n$-character 
string over an alphabet of size $k$, and which is the unique minimum element 
in the lexicographical ordering of all its rotations. 
Being the singularly smallest rotation implies that a Lyndon word 
differs from any of its non-trivial rotations, and is therefore aperiodic 
(see \cite{bp}). For example (see \cite{bp}), the list of Lyndon words 
of length 6 on the alphabet $\{0,1\}$ reads
 $$
000001, 000011, 000101, 000111, 001011, 001111, 001111, 010111, 011111.
$$  
Of course, the number of Lyndon words of length $n$ on $k$ symbols
is equal to $M_k(n)$, where $M_k(n)$ is given by  (1). 

Notice that the authors of the paper \cite{bp} investigate the
historical roots of the field of combinatorics of words.  They comprise 
applications and interpretations in algebra, geometry, and combinatorial 
enumeration. Combinatorics of words is a comparatively new area of discrete
mathematics. It is pointed out in \cite{bp} that the collective volumes 
written under the pseudonym of Lothaire give an account of it 
(Lothaire's first volume \cite{lo1} appeared in 1983 and was reprinted 
with corrections in 1997 \cite{lo2}).  
 
It is also well known (see, e.g., \cite[p. 162]{ri}) that the number of 
$(n,k)$ -  necklaces is given by 
 $$
N_k(n)=\frac{1}{n}\sum_{d\mid n}\varphi(d)k^{n/d}.\eqno(2)
  $$
The formula (2) is called  
 {\it MacMahon's} formula in the book by Graham et al. 
\cite[the formula (4.63), p. 141]{gr}, while in 
Lucas' book \cite[p. 503]{lu}, it is credited to M. le 
colonel Moreau (see the sequence 
A001031 in   \cite{oe} concerning the sequence 
 $\{N_2(n)\}_{n=1}^{\infty}$) which presents the number of binary necklaces). 
A proof of 
(2) given in \cite[pp. 14-141]{gr} is based on a lemma 
presented by P\'{o}lya \cite{po} (see also Lemma in \cite[p. 659]{ri}). 

Let $G$ be a finite group that acts on a set $X$. For each $g\in G$ let 
$X^g$ denote the set of elements in $X$ that are fixed by $g$. Burnside's 
lemma asserts the following formula for the number of orbits, denoted 
$|X/G|$:
  $$
|X/G|=\frac{1}{|G|}\sum_{g\in G}|X^g|.\eqno(3)
  $$
Two elements of $X$ belong to the same ``orbit'' when one can be reached from 
the other by through the action of an element of $G$. For example,
if $X$ is the set of colorings of a cube, and
$G$ is the set of rotations of the cube, then two elements of  $X$ belong to
the same orbit precisely when one is a rotation of the other.
  
In this paper  we focus our attention to the investigation of enumerations of 
some classes of $n$-character strings and   binary   necklaces, i.e., for 
new deductions of expressions for numbers of some binary  type necklaces. 
As noticed above,  
 binary necklaces are necklaces  in two colors with length 
$n$. Observe that the authors of the paper \cite{fk} exhibit 
a correspondence between the binary cycles on length $n$ and the lexicographic
composition of the integer $n$. Furthermore, in \cite{fk} the authors give
an algorithm for generating all necklaces of a specific density. 
   
The paper is organized as follows. In Section 2 we present the main results 
and related notions necessary for their formulations and proofs.
Section 3 is devoted to the auxiliary results and related notions 
and notations. 
In  Section 4 we 
prove Theorems 1, 1' and  2
 concerning the numbers of three classes of $j$-character strings
(closely related to some classes of binary necklaces of Lyndon words).
Furthermore, by using these results, we give   proofs (of Theorems 3 and 4)   
of the well known formulae for  two classes of  binary 
necklaces without use of Burnside's lemma  and P\'olya enumeration theorem
  (Burnside's lemma is also called 
Burnside's counting theorem, the Cauchy-Frobenius lemma or the 
orbit-counting theorem).

 Notice that the formula (15) of Theorem 3
is a special (binary) case of formula (1)  with $k=2$. Similarly, 
the formula (16)   of Theorem 4 is a special (binary) case of formula (2) 
with $k=2$.  As applications, we obtain some interesting 
congruences involving the sums of certain binomial coefficients 
and the function $\mu(n)$ or $\varphi(n)$ (Corollaries 1 and 2).
 In particular, we obtain two Lucas' type congruences (Corollaries 3 and 4;
see, e.g., \cite{m} and  \cite{m2}). 
The all our main results and their consequences  are given in Section 2.
  
 Methods used in all our proofs of auxiliary results presented in Section 3 
are very 
combinatorial in spirit and they involve the applications of elementary 
number theory. 
By using these auxiliary results, in Section 4 
we give proofs of our results (Theorems 1-4 and Corollaries 3 and 4).

\section{The main results}


Throughout this paper we suppose that  
$k\ge 2$ and  $r\ge 1$ are fixed 
 $k\ge 2$ is an arbitrary fixed integer and  $r\ge 1$   
is an integer  such that $r\le k$. 
Here, as always in the sequel, we will denote by $(n,m)$ the greatest common 
divisor  of positive integers $n$ and $m$, and  by $|{\mathcal S}|$ the cardinality of
 a finite set ${\mathcal S}$. Usualy, denote by  $\mu(n)$ and $\varphi(n)$ 
the  M\"{o}bius  function and the Euler totient function, respectively.
For any positive integers   $j\ge 1$ and $i\ge 1$ with $j\le i$,
denote by ${\mathcal A}(j,i)$ the collection  of all subsets of
$R_i:=\{0,1,\ldots,i-1\}$  that contain 
exactly $j$ elements.
For given set $A_j=\{a_1,a_2,\ldots,a_j\}\in {\mathcal A}(j,i)$,
 denote by  $l=f_i (A_j)>0$
the smallest positive integer for which the set 
$A_j$ is equal to the set $A_j+l$ modulo
$i$, where $A_j+l=\{a_1+l,a_2+l,\ldots,a_j+l\}$. Since
$A_j$ is equal to $A_j+i$ modulo $i$, $f_i(A_j)$ exists and  
$f_i (A_j)\le i$.

For an arbitrary positive integer $s$ with $1\le s\le i$, 
denote by ${\mathcal A}_s(j,i)$ the subset of ${\mathcal A}(j,i)$ consisting
of those sets $A_j\in {\mathcal A}(j,i)$ for which $f_i(A_j)=s$.
Notice that the set ${\mathcal A}_s(j,i)$ may be considered 
as a class of binary necklaces  with   length $j$ whose  properties are 
described above, i.e., for which $f_i(A_j)=s$.
Then  $R_2=\{0,1\}$,  ${\mathcal A}(2,3)=\{\{0,1\},\{0,2\},\{1,2\}\}\}$. 
Then assuming $A_2=\{0,2\}$, we have $A_2+1=\{1,0\}$ modulo 3,
$A_2+2=\{2,1\}$ modulo 3 and $A_2+3=\{0,2\}=A_2$ modulo 3,
whence it follows that $f_3 (A_2)=3$. If for example, $i=j=3$ and 
$A_3=\{0,2,4\}$, then since $A_3+1=\{1,0,2\}$ and 
$A_3+2=\{2,4,0\}=A_3$ modulo 3, we find that $f_3 (A_3)=2$.

Our investigations are motivated 
by  the following question:

{\it Under what conditions on 
integers $k$, $n$, $r$,  $a_1,a_2,\ldots,a_r$, the  
set $A_r=\{a_1,a_2,\ldots,a_r\}$ equals the set 
$A_r+n=\{a_1+n,a_2+n,\ldots,a_r+n\}$ modulo $k$}?

The answer to this question is given by Proposition 2 in  Section 2.
Namely, if $n=f_k(A_r)$, then by Proposition 1,  $n$ necessarily  divides
$k$,  that is, $k=nd$, and by Proposition 2, $d$ divides $r$, i.e., 
$r=md$ for a positive 
integer $m$. When this is the case, by Proposition 2, the set $A_r$ 
  contains exactly $m$ distinct representatives modulo
$n$. Further,  assuming that $a_1',a_2',\ldots,a_m'$ are 
these  representatives modulo $n$, then $A_r$ has the form 
 $$
A_r=\{a_j'+sn :\, 1\le j\le m,\, 0\le s\le d-1\}.\eqno(4)
 $$
\indent For given common divisor $d$ of $k$ and $r$, with $k=nd$ and $r=md$,
by Proposition 3, we have  $|{\mathcal A}_n(m,n)|=|{\mathcal A}_n(r,k)|$.
This result is the basic tool for determining the cardinality
$|{\mathcal A}_n(r,k)|$
(proof of Theorem 1 in Section 3).

Now we present our basic  result whose proof will be given in Section 4.

\vspace{2mm}

T{\scriptsize HEOREM} 1.
 {\it Let $A_r=\{a_1,a_2,\ldots,a_r\}\subseteq 
\{0,1,2,\ldots,k-1\}$ be a set for which $n=f_k(A_r)$.
Then $n$ divides $k$ and $d=k/n$ is a positive integer that
 divides $r$, that is, $k=nd$ and 
$r=md$ for a positive integer $m$. Furthermore, the class of 
$r$-character strings ${\mathcal A}_n(r,k)$ consists of  
    $$
|{\mathcal A}_n(r,k)|=\sum_{s|(n,m)}{\frac{n}{s}\choose\frac{m}{s}}\mu(s)
\eqno(5)
 $$
elements, and the sum is taken over all positive divisors $s$ of the greatest common divisor $(n,m)$ 
of $n$ and $m$.

In particular, we have}
$$
|{\mathcal A}_k(r,k)|=\sum_{s|(k,r)}{\frac{k}{s}\choose\frac{r}{s}}\mu(s).
\eqno(6)
 $$
\vspace{2mm}

R{\scriptsize EMARK} 1. Notice that the numbers 
$|{\mathcal A}_n(r,k)|$ are closely related to the sequence 
(triangular array read by rows) A185158 in \cite{oe}, Namely, 
$T(n,m)=|{\mathcal A}_n(r,k)|$ for all $n,k=1,2,\ldots$ 
(with $k=nd$ and $r=md$), where by  Comments in \cite{oe}, $T(n,m)$  is 
the number of binary Lyndon words of length $n$ containing $m$ ones
(cf. the sequence/triangular array  A051168 in \cite{oe}).

\vspace{2mm}

R{\scriptsize EMARK}  2. It follows from Theorem 1 and Example 1  at the end
of Section 3 that the collection ${\mathcal A}_n(r,k)$ is  a nonempty   set 
if and only if $d=k/n$ is an integer that divides $r$.

\vspace{2mm}

For a fixed positive integer $n$ that  divides $k$ such that 
the integer $d=k/n$ divides $r$,  
the $r$-character strings that belong to  ${\mathcal A}_n(r,k)$  can be 
separated into disjoint classes as follows. 
 We say that the $r$-character strings $A$ and $A'$ in ${\mathcal A}_n(r,k)$ 
are $k$-{\it equivalent}, writting $A\sim_k A'$,  
if there exists   $j\in\{0,1,2,\ldots,k-1\}$ such that $A+j$ is equal to $A'$ 
modulo $k$. It is easy to see that  $\sim_k$ is an   equivalence relation, and
that every coset (with respect to this relation) 
has exactly $n$ elements. More precisely, the coset $\overline{A}$ represented
by a set $A\in {\mathcal A}_n(r,k)$ is equal  to 
$\{A+j:\,j=0,1,\ldots,n-1\}$ modulo $k$.
Denote by ${\mathcal P}_n(r,k)$ the set of
all these  cosets. We say  that each element of ${\mathcal P}_n(r,k)$ 
is a $(r,k)$-{\it period with length}  $n$. Thus by (5) and (6) 
of Theorem 1, with $k=nd$ and $r=md$, 
we obtain the following result.

\vspace{2mm}

T{\scriptsize HEOREM} 1'.
{\it Suppose that  $k=nd$ and $r=md$ for some positive integers  
$n,$ $m$ and $d$.
Then we have
  $$
|{\mathcal P}_n(r,k)|=\frac{|{\mathcal A}_n(r,k)|}{n}=\frac{1}{n}
\sum_{s|(n,m)}{\frac{n}{s}\choose\frac{m}{s}}\mu(s),
\eqno(7)
 $$
where  the sum is taken
over all positive divisors $s$ of the greatest common divisor $(n,m)$ of $n$ 
and $m$.

In particular, the number of $(r,k)$-periods with maximal length  $k$
is given as}
$$
|{\mathcal P}_k(r,k)|=\frac{|{\mathcal A}_k(r,k)|}{k}=\frac{1}{k}
\sum_{s|(k,r)}{\frac{k}{s}\choose\frac{r}{s}}\mu(s).
\eqno(8)
 $$    

Now define the sum $S(r,k)$ as  
$$
 S(r,k)=\sum_{n=1}^{k}|{\mathcal P}_n(r,k)|=\sum_{n=1}^{k}
\frac{|{\mathcal A}_n(r,k)|}{n},\eqno(9)
$$
 that is, for a fixed $r$, $S(r,k)$ 
is a  number of all $(r,k)$-periods with arbitrary length $n$ ($1\le n\le k$).
Observe that by Theorem 1, $S(r,k)$ may be written as  
 $$
 S(r,k)=\sum_{n\mid k\atop (k/n)|r}\frac{|{\mathcal A}_n(r,k)|}{n}.\eqno(10)
$$

T{\scriptsize HEOREM} 2.    {\it $S(r,k)$ is given as a sum
     $$
S(r,k)=\frac{1}{k}\sum_{s|(k,r)}{\frac{k}{s}\choose\frac{r}{s}}
\varphi(s),\eqno(11)
 $$
 where the sum 
is taken
 over all positive divisors $s$ of the greatest common divisor $(k,r)$ of $k$ 
and $r$.}

\vspace{2mm}

R{\scriptsize EMARK} 3.  It is known (see, e.g., the sequence 
 $L_2(n,d)$ in 
\cite{http})  that the 
number $S(r,k)$ given by (11) in Theorem 2 is the number of 
 binary Lyndon words of length $k$ containing $r$ ones.

\vspace{2mm}

 Notice that the numbers 
$|{\mathcal A}_n(r,k)|$ are closely related to the sequence 
(triangular array read by rows) A185158 in \cite{oe}, Namely, 
$T(n,m)=|{\mathcal A}_n(r,k)|$ for all $n,k=1,2,\ldots$ 
(with $k=nd$ and $r=md$), where by   Comments in \cite{oe}, $T(n,m)$  is 
the number of binary Lyndon words of length $n$ containing $mk$ ones
(cf. the sequence/triangular array  A051168 in \cite{oe}).

As  the immediate consequences of (8) and (11) 
 we get the following  two congruences,  respectively.  

\vspace{2mm}

C{\scriptsize OROLLARY} 1. {\it Let $k\ge 2$ and $r\ge 1$ be  
integers  with $r\le k$. Then we have}
        $$
\sum_{s|(k,r)}{\frac{k}{s}\choose\frac{r}{s}}\mu(s)\equiv 0\pmod{k}.
\eqno (12)
 $$

C{\scriptsize OROLLARY} 2.  {\it Let $k\ge 2$ and $r\ge 1$ be  
positive integers  with $r\le k$. Then we have}
       $$
\sum_{s|(k,r)}{\frac{k}{s}\choose\frac{r}{s}}\varphi (s)\equiv 0\pmod{k}.
\eqno (13)
 $$
\vspace{1mm}

 Let $R(k)$ be the sum defined as  
 $$
R(k)=\sum_{r=1}^{k}|{\mathcal P}_k(r,k)|\eqno (14)
  $$ 
that is, for a fixed $k$, $R(k)$ 
is a number of all  $(r,k)$-periods $(1\le r\le k)$ with arbitrary length
$n$  ($1\le n\le k$),

The following result is a special (binary) case of 
Moreau's necklace-counting function (1) with $k=2$.

\vspace{2mm}

T{\scriptsize HEOREM} 3.
{\it Let $k\ge 2$ be any integer. Then 
 $$
R(k)=M_2(k)=\frac{1}{k}\sum_{s|k}2^{\frac{k}{s}}\mu(s),\eqno(15),
  $$
where $M_2(k)$ is the binary case of 
Moreau's necklace-counting function  given by $(1)$,
 where the sum on the right hand side ranges over all divisors $s$ of $k$.}

\vspace{2mm}

Finally, put $L(k)=\sum_{r=1}^{k} S(r,k)$, that is, for a fixed $k$, $L(k)$ 
is a  number of all $(r,k)$-periods  $(1\le r\le k)$ 
with any possible  length $n$ ($1\le n\le k$).

The following result is a special (binary) case  
of MacMahon's formula (2) with $k=2$.

\vspace{2mm}

T{\scriptsize HEOREM} 4. {\it  Let $k\ge 2$ be any integer. Then 
   $$
L(k)=N_2(k)=\frac{1}{k}
\sum_{s|k}2^{\frac{k}{s}}\varphi(s),\eqno(16)
$$ 
where $N_2(k)$ is the binary case of MacMahon's formula $(2)$ with $k=2$, and
the sum on the right is taken over all positive divisors $s$ of $k$.}
 \vspace{2mm}

Finally, as consequences of the congruence from Corollary 1, in Section 
4 we prove the following two statements. 

\vspace{2mm}

C{\scriptsize OROLLARY} 3. (cf. \cite{m2}) {\it Let $n\ge 1$ and $m\ge 1$ be  
relatively prime integers  with $m\le n$. Then for any prime $p$ and
integer $\alpha \ge 1$,}
 $$
{np^{\alpha} \choose mp^{\alpha}}\equiv {n \choose m}\pmod{np}.\eqno(17)
$$
\vspace{1mm}

C{\scriptsize OROLLARY} 4. (see, e.g.,  \cite[the congruence (5) on p. 6]{m})  {\it Let $n\ge 1$ and $m\ge 1$ be  
any positive integers  with $m\le n$. Then for any prime $p$ we have}
 $$
{np \choose mp}\equiv {n \choose m}\pmod{p}.\eqno(18)
$$

\section{The Collections ${\mathcal A}_s(j,i)$ and 
Auxiliary Results}


Let $i$ be a fixed integer greater than 1, and consider an alphabet 
consisting of the numbers $0,1,\ldots,i-1$.  
With this alphabet form all possible $k$-letter words $(a_1,a_2,\ldots,a_k)$,
where $k$ is also fixed.  There are evidently $i^k$ such words in all.
For our purposes, notice that the set  $R_i:=\{0,1,\ldots,i-1\}$ is a 
complete residue system modulo $i$. 
For a finite subset $A$ of $\mathbf{N}_0:=\{0,1,2,\ldots\}$, denote by 
$\overline{A}(i)$ 
the (unique) subset of $R_i$ consisting of all 
$l\in R_i$ for which there is a $a_l\in A$ such that $a_l\equiv l  \pmod{i}$.
In other words, $\overline{A}(i)$  is a set of representatives modulo $i$
(chosen from the set $R_i$) of all elements which belong to $A$.
 For given two finite subsets $A$ and $B$ of $\mathbf{N}_0$
we say that $A$ equals $B$  modulo $i$ if $\overline{A}(i)=\overline{B}(i)$. 
In this case, we shall often write $A=B$ modulo $i$.

For any positive integers  $i\ge 1$ and $j\ge 1$  with $j\le i$,
denote by ${\mathcal A}(j,i)$ the collection  of all subsets of
$\{0,1,\ldots,i-1\}$ that contain exactly $j$ elements.
Given set $A_j=\{a_1,a_2,\ldots,a_j\}\in {\mathcal A}(j,i)$,
and any positive integer $t$, put $A_j+t=\{a_1+t,\ldots,a_j+t\}$.
 Denote by  $l=f_i (A_j)>0$ the  smallest positive integer for which 
$A_j=A_j+l$ modulo $i$. Clearly, $A_j=A_j+i$ modulo $i$, 
whence we see that $f_i(A_j)$ exists and $f_i (A_j)\le i$. 

For an arbitrary positive integer $s$ with $1\le s\le i$, 
denote by ${\mathcal A}_s(j,i)$ a subset of ${\mathcal A}(j,i)$ consisting
of those sets $A_j$ in  ${\mathcal A}(j,i)$ for which $f_i(A_j)=s$.
It is of interest here to consider the collections ${\mathcal A}_n(r,k)$
and ${\mathcal A}_n(m,n)$.

Recall that  $k\ge 2$ is any fixed integer and   $r\ge 1$   
is an integer  such that $r\le k$.  In this section we give 
necessary conditions on integers $k,r,n,a_1,a_2,\ldots,a_r$, to be satisfied
$f_k(A_r)=n$ for given set $A_r=\{a_1,a_2,\ldots,a_r\}\in {\mathcal A}(r,k)$.
To solve this problem, we start with the following proposition.

\vspace{2mm}

P{\scriptsize ROPOSITION} 1.
{\it For any set $A_r=\{a_1,a_2,\ldots,a_r\}\subseteq 
\{0,1,2,\ldots,k-1\}$, the integer $n=f_k (A_r)$ divides $k$.}
\vspace{2mm} 

{\it Proof}.
As noticed above, $n=f_k (A_r)\le k$. 
If we suppose that the integer $f_k(A_r)=n$ does not divide  
$k$, then $k=q_1n+r_1$ with positive integers $q_1$ and $r_1$such that
  $0< r_1\le n-1$, and hence
$$
A_r=A_r+k=A_r+(q_1n+r_1)=(A_r+q_1n)+r_1=A_r+r_1\quad{\rm modulo}\,\, k.
$$
It follows that $n=f_k(A_r)\le r_1<n$. 
This contradiction shows that $f_k(A_r)$ divides $n$.
\hfill$\qed$
\vspace{2mm} 

Let $d$ be any  divisor of  $k$, and  $k=nd$ for an integer $n\ge 1$.
For a fixed integer $i$ such that $0\le i\le k-1$, consider the set 
$C_i$ defined as
 $$
C_i=\{i,i+n,i+2n,\ldots, i+(d-1)n\}.\eqno(19)
 $$

Then we have the following lemma.

\vspace{2mm} 

L{\scriptsize EMMA} 1. {\it Every set $\overline{C}_i(k)$ 
$(0\le i\le k-1)$ has exactly $d$ elements.
Moreover, $\overline{C}_i(k)=\overline{C}_j(k)$ if and only if 
$i\equiv j \pmod{n}$. In the case when $i\not\equiv\, j\pmod{n}$,
$\overline{C}_i(k)$ and $\overline{C}_j(k)$ are disjoint sets}.

\vspace{2mm} 
 
{\it Proof}. First observe  that the set $\overline{C}_i(k)$ has
$d$ elements. Namely, if $i+d_1n\equiv i+d_2n \pmod{k}$ for some
$d_1$ and $d_2$ with $0\le d_1<d_2\le d-1$, then $(d_2-d_1)n\equiv 
0\pmod{dn}$. Thus $(d_2-d_1)\equiv 
0\pmod{d}$, and hence it must be $d_1=d_2$.

Suppose that $\overline{C}_i(k)$ and $\overline{C}_j(k)$ have
at least one common element. In other words, assume that
$i+d_1n\equiv j+d_2n\pmod{k}$ for some
$d_1$ and $d_2$ with $0\le d_1,d_2\le d-1$, or equivalently,
$i-j+(d_1-d_2)n\equiv 0\pmod{dn}$. Therefore, we obtain
$i\equiv j\pmod{n}$, whence it follows easily that 
$\overline{C}_i(k)=\overline{C}_j(k)$.
\hfill$\qed$

\vspace{2mm}

 Given set $A_r=\{a_1,a_2,\ldots,a_r\}\subseteq 
\{0,1,2,\ldots,k-1\}$, put $n=f_k(A_r)$.  Then by Proposition 1,
 the number $d=k/n$ is an integer. Assume that the set $A_r$ contains $m$ 
distinct representatives modulo $n$. Choose a maximal subset  
$B_m=\{a'_1,a'_2,\ldots,a'_m\}$ of  $A_r$ such that
$a'_s\not\equiv\, a'_q\pmod{n}$ for any integers $s$ and $q$ with 
$1\le s<q\le m$. Then by  Lemma 1,
    $$
\bigcup_{j=1}^r\overline{C}_{a_j}(k)=\bigcup_{j=1}^m\overline{C}_{a'_j}(k),
\eqno(20)
 $$
where  the sets $\overline{C}_{a'_j}(k)$ are disjoint in pairs, that is, 
 $\overline{C}_{a'_s}(k)\bigcap \overline{C}_{a'_q}(k)$ is the empty set
   for any integers $p$ and $q$ with $1\le s<q\le m$. Furthermore,
 $$
\left|\bigcup_{j=1}^r\overline{C}_{a_j}(k)\right|=\left|\bigcup_{j=1}^m
\overline{C}_{a'_j}(k)
\right|=md,\eqno(21)
 $$
Since by the assumption $A_r+n=A_r$ modulo $k$, and hence 
$A_r+ln=A_r$ modulo $k$ for all integers $l$ with  $0\le l\le d-1$, we have
 $$
\bigcup_{j=1}^m\overline{C}_{a'_j}(k)\subseteq \{a_1,a_2,\ldots,a_r\},
\eqno(22)
 $$
or equivalently,
 $$
A_r':=\{a'_i +jn:\, 1\le i\le m,\, 0\le j\le d-1\}\subseteq A_r.\eqno(23)
 $$

\indent On the other hand, for any $a_i\in A_r$, there exists $a'_j$ with
$1\le j\le m$, such that $n$ divides $a_i- a'_j$. Hence, in view of the fact that
$0\le a_i- a'_j\le k-1$, there is an $s$ with $0\le s\le d-1$ such that
$a_i- a'_j=sn$, i.e., $a_i= a'_j+sn\in A'_r$. Therefore, $A_r\subseteq A'_r$,
and hence, it must be $A_r=A'_r$.  It follows that
$r=|A_r|=|A'_r|=md$, and we have
 $$
A_r=\{a'_i +jn:\, 1\le i\le m,\, 0\le j\le d-1\}.\eqno(24)
 $$
The above arguments together with Proposition 1 imply the 
following result.

\vspace{2mm} 
P{\scriptsize ROPOSITION} 2.
{\it For given set 
$A_r=\{a_1,a_2,\ldots,a_r\}\subseteq \{0,1,2,\ldots,k-1\}$  put $n=f_k(A_r)$.
Then $d=k/n$ is an integer that divides $r$, that is, $r=md$ for a positive 
integer $m$. 
Moreover, 
the set $A_r$ contains exactly $m$ distinct representatives modulo
$n$. If we assume that $a_1',a_2',\ldots,a_m'$ are these representatives 
modulo $n$, then $A_r$ has the form} 
$$
A_r=\{a'_i +jn:\, 1\le i\le m,\, 0\le j\le d-1\}.\eqno(25)
 $$

\vspace{2mm}

R{\scriptsize EMARK} 4. Clearly, $A_r=A_r+n$ modulo $k$
for every set $A_r$ given by (25). However, the converse of Proposition 2 
is not true in the sense that  generally, given positive integers $k,r,n,m$ 
and $d$ such that $k=nd$ and $r=md$, there are sets $A_r$ of the form
(25) for which $f_k(A_r)<n$. To show this fact, put $n=k=4$, $d=1$, $m=r=2$,
 and consider the set $A_2=\{0,2\}\subset \{0,1,2,3\}$.
Then $A_2+2=A_2$ modulo $4$, and hence $f_4(A_2)=2<4$.

\vspace{2mm}
R{\scriptsize EMARK} 5.
If the integers $k$ and $r$ are relatively prime,
using the same notations as in Proposition 2, 
this proposition implies that $d$ divides $(k,r)=1$. Hence, it must be
$d=1$ and $k=n=f_k(A_r)$ for any set $A_r=\{a_1,a_2,\ldots,a_r\}\in 
{\mathcal A}(r,k)$. It follows that  $\left|{\mathcal A}_k(r,k)\right|
=\left|{\mathcal A}(r,k)\right|={k\choose r}$
(cf. (5) and (6) of Theorem 1). This means 
that  each set $A_r\in {\mathcal A}(r,k)$   belongs
to  certain  $(r,k)$-period with maximal length  $k$.

\vspace{2mm}
The following result has an important role in the proof of Theorem 1.

\vspace*{2mm}

P{\scriptsize ROPOSITION} 3.
 {\it
For  an arbitrary common divisor $d\ge 1$
of $k$ and $r$, take $k=nd$ and $r=md$.
Then the collections ${\mathcal A}_n(m,n)$ and ${\mathcal A}_n(r,k)$ have same 
cardinality, and one bijection $h$ between these collections is given as}
   $$
B_m=\{a_1,a_2,\ldots,a_m\}\mapsto A_r=\{a_i+jn:\, 1\le i\le m,\, 
0\le j\le d-1\},\eqno(26)
  $$
{\it where $B_m$ is in ${\mathcal A}_n(m,n)$
and $A_r$ is in ${\mathcal A}_n(r,k)$.}

\vspace{2mm}

{\it Proof}. 
For a given set $B_m=\{a_1,a_2,\ldots,a_m\}\in {\mathcal A}_n(m,n)$,
it is easy to check  that all elements
of its  associated set $A_r=\{a_i+jn:\, 1\le i\le m,\, 
0\le j\le d-1\}$ are distinct modulo $k$. Therefore, $|A_r|=md=r$ modulo $k$,
and hence the above map $B_m\mapsto A_r$, denoted here as $h$, 
 is a map into ${\mathcal A}(r,k)$.
Furthermore,  it is routine to verify that the map $h$  is injective.

It remains to show that the map $h$ is onto ${\mathcal A}_n(r,k)$.
Let $A_r\in {\mathcal A}_n(r,k)$ be  arbitrary. 
This means that $n=f_k(A_r)>0$ is the 
smallest positive integer for
which $A_r+n$ equals $A_r$ modulo $k$.  It follows from Proposition 2  that 
the set $A_r$ contains  $m$ distinct representatives modulo $n$; assume  
$a_{i_1}, a_{i_2},\ldots, a_{i_{m}}$. Hence, by (25) we get
 $$
A_r=\{a_{i_l}+jn :\, 1\le l\le m,\, 0\le j\le d-1\}.
\eqno(27)  
 $$
Obviously, it is sufficient to show that there exists a set 
$B_m\in {\mathcal A}_n(m,n)$ such that $h(B_m)=A_r$. 
Define $B_m=\{a_{i_1},a_{i_2},\ldots,a_{i_m}\}$. 
Clearly, $B_m$ is in ${\mathcal A}(m,n)$,
and hence it suffices to show that $f_n(B_m)=n$. Suppose that 
$f_n(B_m)=n_1<n$. 
Then as in  the proof of Proposition 1,  we infer that $n_1$ divides $n$, 
i.e., $n=n_1d_1$ with an integer $d_1>1$.
Since the set $B_m+n_1$ equals $B_m$ modulo $n$,   
by  Proposition 2, with $n_1$, $m$, $n$,  $d_1$ and $B_m$ instead of 
$n$, $r$ $k$, $d$ and $A_r$, respectively, we conclude that 
$d_1$ divides $m$, i.e., $m=m_1d_1$ with  $m_1\in\mathbf{N}$. Furthermore, 
by Proposition 2, the set $B_m$ contains exactly $m_1$ distinct 
representatives modulo $n_1$,  assume for example, 
$a_{i_1}, a_{i_2},\ldots, a_{i_{m_1}}$. Then by (25)  of 
Proposition 2, $B_m$ has the form
     $$
B_m=\{a_{i_s}+tn_1 :\, 1\le s\le m_1,\, 0\le t\le d_1-1\},\eqno(28)  
 $$
which by (27) implies that 
$$
A_r=\{(a_{i_s}+tn_1)+jn\, :\, 1\le s\le m_1,\, 0\le t\le d_1-1,\,
0\le j\le d-1\},  \eqno(29)
 $$
whence by putting $n=n_1d_1$, we obtain
  $$
A_r=\{a_{i_s}+(t+jd_1)n_1 \,:\, 1\le s\le m_1,\, 0\le t\le d_1-1,\,
0\le j\le d-1\}.  \eqno(30)
 $$
Because of  $k=nd=n_1d_1d$, it is easy by (30) to verify that the 
set $A_r+n_1$ is equal to the set $A_r$  modulo $k$. This implies that 
$f_k(A_r)\le n_1<n$. This contradiction 
with our assumption that $f_k(A_r)=n$ shows that $f_n(B_m)=n$. This means that
$B_m$ is in ${\mathcal A}_n(m,n)$, and since by (26), (27) and the definition 
of $B_m$, $h(B_m)=A_r$, we conclude that $h$ is a surjective map.
This  completes the proof.
\hfill$\qed$

\vspace{2mm}

E{\scriptsize XAMPLE} 1.  We will show that ${\mathcal A}_n(r,k)$ is a  nonempty 
set for each positive integer $n$ satisfying conditions of Proposition
2. More precisely, for any integers $n\ge 1$ and $d\ge 1$, 
such that $k=nd$ and $r=md$, we will construct some elements of 
${\mathcal A}_n(r,k)$. 
Since by (26) of Proposition 3 it is given  an one-to-one
correspondence $h$ between the families ${\mathcal A}_n(m,n)$  
and ${\mathcal A}_n(r,k)={\mathcal A}_n(md,nd)$,
it is sufficient to consider the corresponding problem for the
families ${\mathcal A}_n(m,n)$ with $m\le n$. 

If $(n,m)=1$, then according to  Remark 5
(by replacing $k$ and $r$ with $n$ and $m$,
respectively), we have ${\mathcal A}_n(m,n)={\mathcal A}(m,n)$. 
In other words, each set $A_m=\{a_1,a_2,\ldots,a_m\}\in {\mathcal A}(m,n)$
is in ${\mathcal A}_n(m,n)$, and hence $\left|{\mathcal A}_n(m,n)\right|={n\choose m}$.

Now we suppose that $(n,m)>1$.  
First we observe that if
the set $A_m$ is in ${\mathcal A}_n(m,n)$, then it is easily seen that the  
 set $A_{n-m}=\{0,1,\ldots,n-1\}\setminus A_m$ is in 
${\mathcal A}_{n}(n-m,n)$. Indeed, since $A_m+n$ is equal to $A_m$ modulo $n$,
then 
    \begin{eqnarray*}
A_{n-m}+n &=&(\{0,1,\ldots,n-1\}\setminus A_m)+n=
\{0,1,\ldots,n-1\}\setminus (A_m+n)\\
\qquad&=&\{0,1,\ldots,n-1\}\setminus A_m=A_{n-m}\,\,{\rm modulo}\,\, n.\qquad\qquad (31)
\end{eqnarray*}

Hence, $f_n(A_{n-m})\le n$. If we suppose that $f_n(A_{n-m})=s<n$,
then as above we obtain that  $A_m+s$ is equal to $A_m$ modulo $n$.
This contradiction with the fact that $A_m$ is in ${\mathcal A}_n(m,n)$
implies that $A_{n-m}$ is in ${\mathcal A}_{n}(n-m,n)$.

In view of the above natural correspondence between
the families ${\mathcal A}_n(m,n)$ and ${\mathcal A}_n(n-m,n)$, and the fact that 
$(n,m)>1$, we may suppose that $2\le m\le\lfloor\frac{n}{2}\rfloor$. 
For such a $m$ define
    $$
V=\{v\in{\rm{\bf N}}:\, v\le 2m-2,\,(v,n)>1\}.\eqno(32)
$$
Since $v=m\le 2m-2$ and by the assumption, $(n,m)>1$, 
it follows that $V$ is a  nonempty set. Let $u=\max\{v|\, v\in V\}$ and
 $$
S=\{s_1s_2:\, s_1\in {\rm{\bf N}}, s_2\in {\rm{\bf N}}, s_1<\frac{n}{u},
s_2\le n, (s_2,n)=1\}.\eqno(33)
 $$
We will prove that $B_s=\{s,2s,\ldots, ms\}\in {\mathcal A}_n(m,n)$
for each $s\in S$. In particular, we have 
 $\{1,2,\ldots, m\}\in {\mathcal A}_n(m,n)$.
 First show that $B_s$ has exactly $m$ different elements modulo $n$. 
Assume that for a fixed $s\in S$, the integers $l_1s$ and $l_2s$ are in 
$B_s$ with 
$1\le l_1-l_2\le m-1$ such that $n\mid (l_1-l_2)s$. We write $s=s_1s_2$ with 
integers $s_1$ and $s_2$ as described by  (33). Then 
since  $(s_2,n)=1$, it must be $n\mid (l_1-l_2)s_1$. Now,
since $1\le s_1<\frac{n}{u}< n$, it follows that  $l_1-l_2=l_3s_3$ 
such that $(s_3,n)=1$ and $(l_3,n)>1$. So $n\mid l_3s_1$, and because of 
$l_3\le m-1$ and $(l_3,n)>1$, we have $l_3\le u$ with $u$  defined above. 
Thus $1\le l_3s_1<u\cdot \frac{n}{u}=n$, and therefore, 
$l_3s_1 \not\equiv\, 0{(\bmod\, n)}$.
This contradiction shows that $|S|=m$ modulo $n$.

It remains to show that $B_s$ is not equal to
$B_s+t$ modulo $n$ for any integer $t$ with $1\le t\le n-1$.
Indeed, if  $B_s=B_s+t$ modulo $n$ for some $t$ with $1\le t\le n-1$,
then there is an integer $p$ with $1\le p\le m-1$ such that $s+t=(p+1)s$.
Therefore, $B_s+t= \{(p+1)s, (p+2)s,\ldots,(p+m)s\}$ modulo $n$, and hence 
we have
   $$
\{s,2s,\ldots, ms\}=\{(p+1)s, (p+2)s,\ldots,(p+m)s\}\,\,{\rm modulo}\,\, n.
\eqno(34) 
$$
Thus there is an integer $q$ with $p+1\le q\le p+m\le 2m-1$ such that
$n\mid qs-s$. It follows that 
$n\mid (q-1)s_1s_2$, whence
since $(s_2,n)=1$, we have $n\mid (q-1)s_1$. Since $s_1\le n-1$, 
it follows that $(q-1,n)>1$, i.e., $q-1=q_1q_2$ for integers $q_1$ and $q_2$ 
with $1\le q_1\le q-1\le 2m-2$ such that $(q_1,n)>1$ and $(q_2,n)=1$. 
This yields that $q_1\le u$, and $n\mid  q_1s_1$, which is impossible since
$1\le q_1s_1<u\cdot\frac{n}{u}=n$. This contradiction implies 
that  $B_s$ is in ${\mathcal A}_n(m,n)$ for each $s\in S$.

\section{Proofs of Theorems 1--4  and Corollaries 3 and 4}


We give here a combinatorial proof of Theorem 1 which is 
based on auxiliary results obtained in Section 3 and on property of function 
$\tau (a,t)$ defined as follows.

\vspace{2mm}

D{\scriptsize EFINITION} 1. For any integers $a>1$ and $t\ge 1$, denote
by $\tau (a,t)$ the number of $t$-tuples $(a_1,a_2,\ldots,a_t)$ of integers
such that  $ a_j\ge 2$ for all $j=1,2,\ldots,t$, and $a=a_1a_2\cdots a_t$. 
Obviously,   $\tau (a,t)=0$ for all $t\ge a$.
\vspace{2mm} 

In the proof of Theorem 1, we use  the following 
property of the function $\tau (a,t)$. 
\vspace{2mm}

L{\scriptsize EMMA} 2. {\it For each  integer $a>1$, we have
     $$
\sum_{t=1}^{a-1}(-1)^t\tau (a,t)=\mu(a), \eqno(35)
    $$
where $\mu(a)$ is the M\"{o}bius function. }

\vspace{2mm} 

{\it Proof}. We derive the proof by induction on $a\ge 2$.
Since $-\tau (2,1)=-1=\mu(2)$, we see that (35) is true for $a=2$. 
Suppose that $a>2$ and that (35) is satisfied for all integers less than 
$a$. 

Obviously, there holds $\tau (a,1)=1$ for all $a>1$. Letting that  the first 
coordinate $q=a_1$ of $t$-tuples $(a_1,a_2,\ldots,a_t)$ of integers
satisfying $a_j\ge 2$ for all $j=1,2,\ldots,t$, and $a=a_1a_2\cdots a_t$,
is taken
over all divisors of $a$, 
by Definition 1 of  $\tau (a,t)$, we have
     $$
\tau (a,t)=\sum_{q|a\atop q>1}\tau\left(\frac{a}{q},t-1\right).\eqno(36)   
  $$
Now  by using the induction hypothesis, (36) and the basic property of 
the M\"{o}bius function (see, e.g., [1, (32) on p. 181]) given by
  $$
\sum_{s|a}\mu(s)=\left\{\begin{array}{rl}
    1 & \qquad {\rm if}\,\, a=1 \\
    0 & \qquad {\rm if}\,\, a>1,
           \end{array}\right.\eqno(37) 
      $$
 we get
    \begin{eqnarray*}
\qquad\qquad \sum_{t=1}^{a-1}(-1)^t\tau (a,t)&=&\sum_{t=2}^{a-1}(-1)^t
\sum_{q|a\atop 1<q<a}\tau\left(\frac{a}{q},t-1\right)-\tau(a,1)\\
&=&-\sum_{t=2}^{a-1}
\sum_{q|a\atop 1<q<a}(-1)^{t-1}\tau\left(\frac{a}{q},t-1\right)-1\\
\qquad\qquad\qquad &=&-\sum_{q|a\atop 1<q<a}\sum_{t=1}^{a-1}(-1)^{t}\tau\left(\frac{a}{q},t\right)
-1\qquad\qquad\qquad\quad (38)\\
&=&\qquad\qquad
-\sum_{ q|a\atop 1<q<a}\mu\left(\frac{a}{q}\right)-1
\,\, ({\rm by\,\, the\,\, hypothesis\,\, for\,\,}\frac{a}{q}<a)\\
&=&-\sum_{q|a\atop 1\le q\le a}\mu\left(\frac{a}{q}\right)+\mu(a)
= -\sum_{s|a}\mu(s)+\mu(a)=\mu(a).
    \end{eqnarray*}
This completes the proof. 
\hfill$\qed$
\vspace{2mm}

We are now ready to prove the main result.

\vspace{2mm}

{\it Proof of Theorem 1}. Note that the first assertion 
of Theorem 1 is contained in Proposition 2. 
It remains to prove the equality (5).

Using the notations introduced in Section 3, if 
$k=nd$ and $r=md$ for an integer $d\ge 1$,  then by Proposition 2 we have 
     $$
|{\mathcal A}_n(md,nd)|=|{\mathcal A}_n(m,n)|.\eqno(39)
     $$
If  integers $k$ and $r$ are relatively prime, 
then it must be $d=1$, $r=m$, and thus $n=k=f_k(A_r)$ for any set 
$A_r=\{a_1,a_2,\ldots,a_r\}\subset \{0,1,2,\ldots,k-1\}$. Therefore,  
$$
\left|{\mathcal A}_n(r,k)\right|
=\left|{\mathcal A}(r,k)\right|={k\choose r}={n\choose m}=\mu(1)=
{n\choose m},\eqno(40)
 $$
whence follows (5). 

Now suppose that $(k,r)>1$. To determine $|{\mathcal A}_n(r,k)|$,
where $k=nd$ and $r=md$ with $d\ge 1$, denote
   $$
\overline{{\mathcal A}}_n(m,n)={\mathcal A}(m,n)\setminus {\mathcal A}_n(m,n).
\eqno(41)
   $$
Then since $|{\mathcal A}(m,n)|={n \choose m}$, we have
   $$
|{\mathcal A}_n(m,n)|={n \choose m} -|\overline{{\mathcal A}}_n(m,n)|.
\eqno(42)
   $$
Moreover, $\overline{{\mathcal A}}_n(m,n)=\bigcup_{1\le s<n}{\mathcal A}_s(m,n)$,
and by Proposition 2, ${\mathcal A}_s(m,n)$ is a  nonempty set
 if and only if $s$ divides $n$ and
$d_1=\frac{n}{s}$ divides $m$. When this is the case, by (39),
with $\frac{n}{d_1}$, $\frac{m}{d_1}$ and $d_1$ instead of
$n$, $m$ and $d$, respectively, we obtain
   $$
|{\mathcal A}_{\frac{n}{d_1}}(m,n)|=\left|{\mathcal A}_{\frac{n}{d_1}}\left(\frac{m}{d_1},
\frac{n}{d_1}\right)\right|.\eqno(43)
 $$
Therefore, we obtain
 $$
\overline{{\mathcal A}}_n(m,n)=\sum_{d_1|(m,n)
\atop d_1>1}\bigcup{\mathcal A}_{\frac{n}{d_1}}
\left(\frac{m}{d_1},\frac{n}{d_1}\right),\eqno(44) 
  $$
whence it follows that
    $$
|\overline{{\mathcal A}}_n(m,n)|=\sum_{d_1|(m,n)
\atop d_1>1}\left|{\mathcal A}_{\frac{n}{d_1}}
\left(\frac{m}{d_1},\frac{n}{d_1}\right)\right|.\eqno(45)
 $$
Thus by  (42) and (45), we get
$$
|{\mathcal A}_n(m,n)|={n \choose m}-\sum_{d_1|(m,n)\atop d_1>1}\left|
{\mathcal A}_{\frac{n}{d_1}}\left(\frac{m}{d_1},\frac{n}{d_1}\right)\right|.
\eqno(46)
 $$
Applying (46) on the all terms of the sum on the right hand side of (46),
 with $\frac{m}{d_1}$ instead 
of $m$ and $\frac{n}{d_1}$ instead of $n$, 
and iterating the same procedure at most $m-1$ times, we have
  \begin{eqnarray*} 
|{\mathcal A}_n(m,n)| &=&{n \choose m}-\sum_{d_1|(m,n)\atop d_1>1}\left|
{\mathcal A}_{\frac{n}{d_1}}\left(\frac{m}{d_1},\frac{n}{d_1}\right)\right|\\
&=&{n \choose m}- \sum_{d_1|(m,n)}{\frac{n}{d_1}\choose\frac{m}{d_1}}\\
&+&\sum_{d_1|(m,n)\atop d_1>1}\sum_{d_2|\left(\frac{m}{d_1},
\frac{n}{d_1}\right)\atop d_2>1}\left({\frac{n}{d_1d_2}\choose\frac{m}{d_1d_2}}
-\left|{{\mathcal A}}_{\frac{n}{d_1d_2}}\left(\frac{m}{d_1d_2},
\frac{n}{d_1d_2}\right)\right|\right)\\
 &=&{n \choose m}- \sum_{d_1|(m,n)\atop d_1>1}{\frac{n}{d_1}
\choose\frac{m}{d_1}}
+\sum_{d_1|(m,n)\atop d_1>1}\sum_{d_2|\left(\frac{m}{d_1},
\frac{n}{d_1}\right)\atop d_2>1}
 {\frac{n}{d_1d_2}\choose\frac{m}{d_1d_2}}\qquad\quad(47)\\
&-&\sum_{d_1|(m,n)\atop d_1>1}
\sum_{d_2|\left(\frac{m}{d_1},\frac{n}{d_1}\right)\atop d_2>1}
\left|{\mathcal A}_{\frac{n}{d_1d_2}}\left(\frac{m}{d_1d_2},\frac{n}{d_1d_2}
\right)\right|\\
&=&\ldots \\
&=&
{n \choose m}\\
&+&\sum_{j=1}^{m-1}(-1)^{j}\left(\sum_{d_1|(m,n)\atop d_1>1}
\sum_{d_2|\left(\frac{m}{d_1},\frac{n}{d_1}\right)\atop d_2>1}\cdots
\sum_{d_j|\left(\frac{m}{d_1\cdots d_{j-1}},\frac{n}{d_1\cdots 
d_{j-1}}\right)\atop d_j>1}
{\frac{n}{d_1d_2\cdots d_{j}}\choose\frac{m}{d_1d_2\cdots d_{j}}}\right).
\end{eqnarray*} 
Hence for a fixed divisor $s>1$  of $(n,m)$ with $s=d_1d_2\cdots d_{j}$ for 
some $j\ge 1$ and the integers $d_1>1,d_2>1,\ldots ,d_j>1$, 
the factor premultiplying the binomial coefficient ${\frac{n}{s}\choose
\frac{m}{s}}$
in the last sum of (47)  is equal to
    $$
\sum_{j=1}^{m-1}(-1)^{j}\tau(s,j),\eqno(48)
 $$
which  is by Lemma 2 equal to  $\mu(s)$. Therefore, by  (47), we obtain
  $$
|{\mathcal A}_n(m,n)|= \sum_{s|(m,n)}{\frac{n}{s}\choose\frac{m}{s}}\mu(s).
\eqno(49)
 $$
This by (39) implies (5), and this completes the proof of Theorem 1.
\hfill$\qed$

\vspace{2mm}

 We will need the following result for the proof of Theorem 2.  
\vspace{2mm}

L{\scriptsize EMMA} 3. [1, (9) on p. 240]
{\it For any integer $q> 1$, we have
    $$
\sum_{s|q}\frac{\mu(s)}{s}=\frac{\varphi(q)}{q},\eqno(50)
 $$
where $\varphi(q)$ is the Euler  totient function.}

\vspace{2mm}

{\it Proof of Theorem 2.} By (7) of Theorem 1', Theorem 1 and (50) of 
Lemma 3, we have
  \begin{eqnarray*}
\qquad S(r,k)&=&\sum_{n=1}^{k}|{\mathcal P}_n(r,k)|=\sum_{d|(k,r)}
\left|{\mathcal P}_{\frac{n}{d}}(r,k)\right|\\
\qquad\qquad  &=&\sum_{d|(k,r)}
\frac{d\left|{\mathcal A}_{\frac{n}{d}}(r,k)\right|}{k}=\frac{1}{k}\sum_{d|(k,r)}d
\sum_{s|\left(\frac{k}{d},\frac{r}{d}\right)}{\frac{k}{ds}\choose
\frac{r}{ds}}\mu(s)\qquad\qquad\qquad (51)\\
&=&\frac{1}{k}\sum_{q|(k,r)}q{\frac{k}{q}\choose
\frac{r}{q}}\sum_{s|q}\frac{\mu(s)}{s}=\frac{1}{k}\sum_{q|(k,r)}{\frac{k}{q}
\choose\frac{r}{q}}\varphi(q),
  \end{eqnarray*} 
as required. \hfill$\qed$

\vspace{2mm}

{\it Proof of Theorem 3}. To determine 
$R(k)=\sum_{r=1}^{k}|{\mathcal P}_k(r,k)|$, by (8) of Theorem 1',
and using the well    known property $\sum\mu(d)_{d\mid n}=1$ if $n=1$, and 
$\sum\mu(d)_{d\mid n}=0$ if $n>1$ (see, e.g., \cite{sh}), 
for each $k>1$ we find that  
    \begin{eqnarray*}
R(k)&=&\frac{1}{k}\sum_{r=1}^{k}\sum_{s|(k,r)}{\frac{k}{s}\choose
\frac{r}{s}}\mu(s)\\
&=&\frac{1}{k}\sum_{s|k}\left({\frac{k}{s}\choose\frac{s}{s}}+
{\frac{k}{s}\choose\frac{2s}{s}}+\ldots+{\frac{k}{s}\choose\frac{k}{s}}\right)
\mu(s)\\
&=&\frac{1}{k}\sum_{s|k}\sum_{i=1}^{\frac{k}{s}}{\frac{k}{s}\choose i}\mu(s)\\
&=&\frac{1}{k}\sum_{s|k}\left(2^{\frac{k}{s}}-1\right)\mu(s)
\qquad\qquad\qquad\qquad\qquad (52)\\
&=&\frac{1}{k}\sum_{s|k}2^{\frac{k}{s}}- \frac{1}{k}\sum_{s|k}\mu(s)
\\
&=&\frac{1}{k}\sum_{s|k}2^{\frac{k}{s}},
  \end{eqnarray*}
as desired. \hfill$\qed$

\vspace{2mm} 

{\it Proof of Theorem 4.} The proof follows in the same manner as that of
Theorem 3 with $\varphi(s)$ instead of $\mu(s)$, by using the well known 
property $\sum_{d\mid n}\varphi(d)=n$ established by Gauss 
(see, e.g., \cite{sh}), and hence may be omitted.
\hfill$\qed$

\vspace{2mm}

{\it Proof of Corollary 3.}  We proceed by induction on $\alpha\ge 1$. 
If $\alpha =1$, then since $(n,m)=1$, $\mu(1)=1$ and $\mu(p)=-1$,  
 (12) of Corollary 1 with $np$ and $mp$ 
instead of $k$ and $r$, respectively, immediately implies that
 $$
{np \choose mp}- {n \choose m}\equiv 0\pmod{np}.\eqno(53)
$$
Now suppose that $\alpha\ge 2$ and (17) holds for all 
positive integers $\beta<\alpha$. Then by using 
the fact that $\mu(p^{\beta})=0$ for each $\beta >1$, (12) gives  
$$
{np^{\alpha} \choose mp^{\alpha}}\equiv {np^{{\alpha}-1} 
\choose mp^{{\alpha}-1}} \pmod{np^{\alpha}}.\eqno(54)
$$
The above congruence together with the induction hypothesis
$
{np^{{\alpha}-1} \choose mp^{{\alpha}-1}}\equiv {n \choose m}\pmod{np}
$
yields  (17). This completes the induction proof. 
\hfill$\qed$
\vspace{2mm}

{\it Proof of Corollary 4.} 
We deduce  the proof by induction on $\sigma =n+m\ge 2$. 
If $\sigma =2$, that is $n=m=1$, (18) is obvious. Suppose that 
$\sigma>2$ and that the congruence (18) is satisfied for
any $n$ and $m$ such that $n+m<\sigma$.  

Assume that $n'$ and $m'$ be positive integers
such that $n'+m'=\sigma$.
If $n'$ and $m'$ are relatively prime, then
(18) is in fact (17) of Corollary 1  with $n=n'p$, $m=m'p$ and $\alpha=1$.  
Now suppose that $(n',m')=d>1$. If $d=p^{\alpha}$ with $\alpha \ge 1$,
i.e., $n'=n''p^{\alpha}$ and $m'=m''p^{\alpha}$ with $(n'',m'')=1$,
then (12) implies that 
$$
{n'p \choose m'p}\equiv{n''p^{\alpha+1} \choose m''p^{\alpha+1}}
\equiv {n'' \choose m''}\pmod{p}.\eqno(55)
$$
If there exists a prime $q\not= p$ that divides $(n',m')=d$, 
then applying the induction hypothesis on integers $n''=n'/q$ and $m''=m'/q$, 
for any divisor $s$ of $\left(n''/q,m''/q\right)$,   we get
 $$
{\frac{n''p}{q}\choose\frac{m''p}{q}}
\equiv {\frac{n''}{q}\choose\frac{m''}{q}}\pmod{p}.\eqno(56)
 $$
By (12) of Corollary 1, we have 
  $$
\sum_{s|(n'p,m'p)}{\frac{n'p}{s}\choose\frac{m'p}{s}}\mu(s)\equiv 0\pmod{p}.
\eqno(57)
 $$
Since $\mu (s'p)=0$ if $p\mid s'$, and  each divisor $s'$ of $(n',m')$ 
with $s'\not\equiv\, 0 {(\bmod\, p)}$ can be uniquely  associated to the divisor 
$s'p$ of $(n'p,m'p)$ with $\mu(s'p)=-\mu(s')$, the above congruence can be written 
as
   $$
{n'p \choose m'p}-{n' \choose m'}+
\sum_{1<s'|(n',m')\atop s'\not\equiv\, 0 {(\bmod\, n)}}
\left({\frac{n'p}{s'}\choose\frac{m'p}{s'}}
-{\frac{n'}{s'}\choose\frac{m'}{s'}}\right)\mu(s')\equiv 0\pmod{p}.\eqno(58)
   $$
Since each term into parantheses is by the  hypothesis divisible by 
$p$, we obtain 
   $$
{n'p \choose m'p}\equiv {n' \choose m'}
\pmod{p}.\eqno(59)
 $$
  This finishes the induction proof. 
\hfill$\qed$

  \end{document}